\documentclass[11pt]{amsart}

\usepackage{amsmath, amssymb, amsthm, amsrefs}
\usepackage{latexsym}
\usepackage{indentfirst}
\usepackage{graphicx}
\usepackage{placeins}
\usepackage{booktabs}
\usepackage{algorithm}
\usepackage{algorithmic}
\usepackage{multirow}

\theoremstyle{plain}
\newtheorem{theorem}{Theorem}

\theoremstyle{definition}

\theoremstyle{remark}

\DeclareMathOperator{\diag}{diag}

\newcommand{\eps}{\epsilon}

\newcommand{\bbm}{\begin{bmatrix}}
\newcommand{\ebm}{\end{bmatrix}}

\newcommand{\p}{\partial}
\newcommand{\G}{\mathcal{G}}
\newcommand{\Q}{\mathcal{Q}}
\newcommand{\F}{\mathcal{F}}

\begin{document}

\title[Directional Computation of High Frequency Boundary
  Integrals]{Fast Directional Computation of High Frequency Boundary
  Integrals via Local FFTs}

\author{Lexing Ying} 

\address{
  Department of Mathematics and Institute for Computational and Mathematical Engineering,
  Stanford University,
  Stanford, CA 94305
}

\email{lexing@math.stanford.edu}

\thanks{This work was partially supported by the National Science
  Foundation under award DMS-0846501 and the U.S. Department of
  Energy’s Advanced Scientific Computing Research program under award
  DE-FC02-13ER26134/DE-SC0009409. The author thanks Anil Damle for
  comments and suggestions.}



\begin{abstract}
  The boundary integral method is an efficient approach for solving
  time-harmonic acoustic obstacle scattering problems. The main
  computational task is the evaluation of an oscillatory boundary
  integral at each discretization point of the boundary. This paper
  presents a new fast algorithm for this task in two dimensions. This
  algorithm is built on top of directional low-rank approximations of
  the scattering kernel and uses oscillatory Chebyshev interpolation
  and local FFTs to achieve quasi-linear complexity. The algorithm is
  simple, fast, and kernel-independent. Numerical results are provided
  to demonstrate the effectiveness of the proposed algorithm.
\end{abstract}


\maketitle

\section{Introduction}

This paper is concerned with the solution of time-harmonic acoustic
obstacle scattering problems. Let $\omega$ be the frequency of the
wave field and $\Omega$ be the scatterer with boundary $\p\Omega$. The
boundary integral method is an appealing approach for this problem
since it typically requires fewer unknowns and provides better
accuracy when compared to the volumetric type methods. The main
computational task of the boundary integral method is the evaluation
of the integral
\begin{equation}
  u(x) = \int_{\p\Omega} G(x,y) f(y) dy,\quad x\in\p\Omega,
  \label{eq:uxcon}
\end{equation}
where $G(x,y)$ is the Green's function of the Helmholtz operator at
frequency $\omega$, or a similar integral where the kernel is replaced
with one of the derivatives of the Green's function.

The numerical treatment of $u(x)$ requires boundary
discretization. For a typical scattering problem, at least a couple of
points per wavelength $\lambda = 2\pi/\omega$ is required in order for
the solution method to achieve reasonable accuracy. If we assume that
both the diameter and the boundary length of $\Omega$ are of order
$\Theta(1)$, then the boundary $\p\Omega$ shall be discretized with a
set $P$ of $n=\Theta(1/\lambda) = \Theta(\omega)$ points. Discretizing
\eqref{eq:uxcon} with an appropriate numerical quadrature scheme
results in the following discrete summation problem: for each $x\in
P$, evaluate
\begin{equation}
  u(x) = \sum_{y\in P} G(x,y) f(y).
  \label{eq:uxdis}
\end{equation}
Direct evaluation of \eqref{eq:uxdis} takes $\Theta(n^2)$ steps. This
can be computationally intensive for high frequency problems where
$\omega$ and $n$ are large. Hence, there is a practical need for
algorithms that can evaluate \eqref{eq:uxdis} rapidly and accurately.

In the past thirty years, a lot of research has been devoted to this
problem and several successful algorithms have been proposed. The
first, and probably the most well-known approach is the high frequency
fast multipole method initiated by Rokhlin in
\cites{rokhlin-1990-rsiest,rokhlin-1993-dfto} and further developed in
\cites{Coifman-1993,darve-2000-fmm,Epton-1995,song-1995-mfma,Song-1997}. This
method follows the structure of the classical fast multipole method
for the Laplace equation and uses fast spherical harmonic transforms
and diagonal forms to construct the multipole and local expansions
efficiently. This method has an $\Theta(n\log n)$ complexity but is
rather technical.

The second method is the butterfly
algorithm proposed by Michielssen and Boag in
\cite{michielssen-1996-mmda}, where the main observation is that the
interaction between two regions under certain geometric configuration
is numerically low-rank even for high frequency scattering
problems. The complexity of this algorithm is $\Theta(n \log^2n)$ with
a relatively large prefactor. 

A third method is the directional FMM-type algorithm proposed by
Engquist and Ying in
\cite{engquist-2007-fdmaoc,engquist-2009-fdahfas}, where the concept
of a directional parabolic separation condition is introduced to
construct directional low-rank approximations between certain square
regions and wedges with narrow opening angle. This algorithm has a
$\Theta(n\log n)$ complexity and the advantage of being
kernel-independent. Other algorithms that rely on this directional
low-rank idea include \cite{bebendorf-2012,messner-2012-fdms}.

While all the above methods are based on hierarchical partitioning of
the boundary, methods based on the fast Fourier transform
\cite{Bleszynski-1996,bruno-2001-afho} are also used widely. These
methods are relatively easy to implement, but have super-linear
complexity like $\Theta(n^\alpha\log n)$ with $\alpha=1.5$ for
example.

In this paper, we introduce a new algorithm for 2D problems. This
algorithm leverages the directional idea of
\cite{engquist-2007-fdmaoc,engquist-2009-fdahfas} and uses oscillatory
Chebyshev interpolation and local FFTs to speed up the
computation. This new method has a $\Theta(n\log^2 n)$ complexity and
is conceptually much simpler than the previous hierarchical
methods. The rest of this paper is organized as follows. Section 2
explained the algorithm in detail. Section 3 provides numerical
results and discussions are given in Section 4.

\section{Algorithm}

\subsection{Notations}

In 2D, the Green's function of the Helmholtz equation is given by
\[
G(x,y) = \frac{i}{4} H^1_0(\omega |x-y|)
\]
and $\lambda = 2\pi/\omega$ is the wavelength. To simplify the
problem, we assume that the domain $\Omega$ has a $C^2$ boundary and
that both its diameter and its boundary length $L$ are of order
$\Theta(1)$ and let $\rho:\partial\Omega \rightarrow [0,L]$ be the
arclength parameterization of the boundary. The extension to
piecewise $C^2$ boundary is quite straightforward and requires little
extra effort.  

We also assume without loss of generality that the length $L$ of
$\p\Omega$ is equal to $4^q \lambda$, where $q>0$ is an integer. The
actual number $4^q$ is not important but it makes the presentation of
the algorithm simpler. Under these assumptions, it is clear that
$\lambda = \Theta(1/4^q)$ and $\omega=\Theta(4^q)$.

Following the typical practice of using a constant number $p$ of
points per wavelength, we discretize the boundary $\p\Omega$ with
a set $P$ of $n=4^q p$ equally-spaced sampling points.

\subsection{Existence of low-rank approximation}

A function $H(x,y)$ with $x\in X$ and $y\in Y$ is said to have a
separated approximation of $n$ terms with relative accuracy $\eps$, if
there exist functions $\{\alpha_i(x)\}_{1\le i\le n}$ and
$\{\beta_i(y)\}_{1\le i\le n}$ such that
\[
|H(x,y) - \sum_{i=1}^n \alpha_i(x) \beta_i(y) | \le \eps |H(x,y)|.
\]

The following result is the theoretical foundation of new algorithm.
\begin{theorem}
Consider two rectangles 
\begin{align*}
  R_1 &= \left[\frac{d}{2}\lambda, \frac{3d}{2}\lambda\right] \times \left[-\frac{h_1}{2}\lambda,\frac{h_1}{2}\lambda\right],\\
  R_2 &= \left[-\frac{3d}{2}\lambda, -\frac{d}{2}\lambda\right] \times \left[-\frac{h_2}{2}\lambda,\frac{h_2}{2}\lambda\right].
\end{align*}
Suppose that $d>1$ and $d\ge 2 \max(h_1,h_2)^2$. Then, for each
$\eps>0$ sufficiently small, there exists an $(\omega,d)$-independent
constant $C(\eps)$ such that $G(x,y) =\frac{i}{4} H^1_0(\omega|x-y|)$
with $x\in R_1$ and $y\in R_2$ has a separated approximation of
$C(\eps)$ terms with relative accuracy $\eps$.
\label{thm:main}
\end{theorem}

\begin{figure}[h!]
  \begin{center}
    \includegraphics[height=1in]{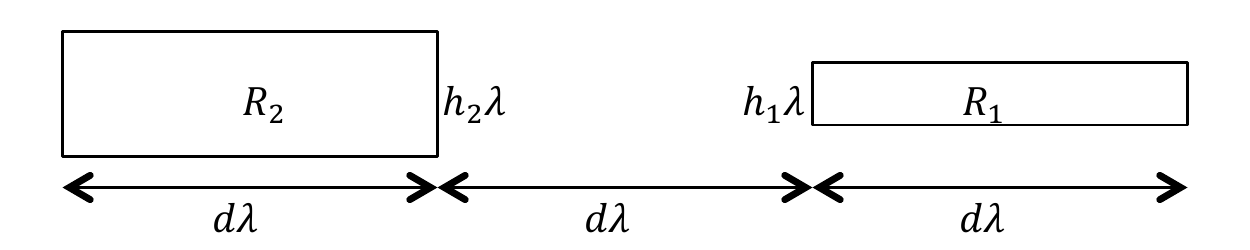}    
  \end{center}
  \caption{The geometric setup of Theorem \ref{thm:main}.}
  \label{fig:thm}
\end{figure}

The main point of the theorem is that the number of terms $C(\eps)$ is
independent of $\omega$ and $d$. Therefore, even if $\omega$ and $d$
grow, the number of terms remains uniformly bounded. For this reason,
such a separated approximation is called low-rank.

\newcommand{\lrsa}{\text{lrsa}}

\begin{proof}
The geometric setting is illustrated in Figure \ref{fig:thm}. We first
observe that $G(x,y)$ can be written as
\[
G(x,y) = \exp(i\omega |x-y|) G_0(|x-y|)
\]
where $G_0(|x-y|)$ is a non-oscillatory real analytic function (see
\cite{abramowitz-1992} for example). A low-rank separation
approximation for $G_0(|x-y|)$ with relative error $\eps$ can be
constructed via Taylor expansion or Chebyshev interpolation. Hence,
the remaining task is to find a low-rank separated approximation for
$\exp(i\omega |x-y|)$ with relative accuracy $\eps$. For $x=(x_1,x_2)
\in R_1$ and $y=(y_1,y_2) \in R_2$, we write
\[
\exp(i\omega|x-y|) = \exp(i\omega(x_1-y_1)) \exp(i\omega(|x-y|-(x_1-y_1))).
\]
As the first term on the right hand side is obviously separated, we
only focus on the second term, which can be rewritten as
\[
\exp(i\omega(|x-y|-(x_1-y_1))) = \exp\left(i\omega(x_1-y_1)
\left(\sqrt{1+\frac{|x_2-y_2|^2}{|x_1-y_1|^2}}-1\right)\right).
\]
We denote the last term by $\eta$ and decompose its calculation into
multiple steps by introducing
\begin{align*}
  & \xi = \frac{1}{|x_1-y_1|^2},\quad
  \alpha = |x_2-y_2|^2\xi,\quad
  \beta = \sqrt{1+\alpha}-1,\quad\\
  & \gamma= i\omega(x_1-y_1)\beta, \quad
  \eta = \exp(i\gamma).
\end{align*}
For a quantity $\delta$, we use the notation $\lrsa(\delta)$ to stand
for a separated approximation of $\delta$ with a number of terms that is
$(\omega,d)$-independent and with a relative error bounded by a
constant multiple of $\eps$.

For $\xi = 1/{|x_1-y_1|^2}$, applying a truncated Taylor expansion of
$x_1$ at the center of $R_1$ with relative error $\eps$ gives a
low-rank separated approximation $\lrsa(\xi)$ for which the number of
terms is $(\omega,d)$-independent due to scale invariance:
\[
\lrsa(\xi) = \xi(1+\eps_\xi), \quad |\eps_\xi|\lesssim \eps.
\]

For $\alpha=|x_2-y_2|^2\xi$, we define $\lrsa(\alpha) = |x_2-y_2|^2
\lrsa(\xi)$. Clearly, by letting $\eps_\alpha = \eps_\xi$, we have
\[
\lrsa(\alpha) = \alpha(1+\eps_\alpha), \quad |\eps_\alpha|\lesssim \eps.
\]
From the relation $\alpha = |x_2-y_2|^2/|x_1-y_1|^2$ and the geometric
setup in Figure \ref{fig:thm}, it is clear that
$\alpha\in(0,1/4)$. Hence, by choosing $\eps$ sufficiently small, one
can easily guarantee that $\lrsa(\alpha)\in(0,1/3)$.

Next, we consider $\beta=\sqrt{1+\alpha}-1$. Since both $\alpha$ and
$\lrsa(\alpha)$ is between $0$ and $1/3$, expanding
$\sqrt{1+\alpha}-1$ with truncated Taylor series at $\alpha=0$ with
relative error $\eps$ gives a polynomial $p(\alpha)$ for which the
number of terms is independent of $\omega$ and $d$. We then define
$\lrsa(\beta) = p(\lrsa(\alpha))$. Notice that
\begin{align*}
  \lrsa(\beta) &= (\sqrt{1+\lrsa(\alpha)}-1)(1+\eps) \\
  &= (\sqrt{1+\alpha(1+\eps_\alpha)}-1)(1+\eps) \\
  &= (\sqrt{1+\alpha}-1)(1+C\eps_\alpha)(1+\eps) \\
  &= (\sqrt{1+\alpha}-1)(1+C\eps_\alpha)(1+\eps)\\
  &= \beta (1+C\eps_\alpha)(1+\eps),
\end{align*}
where $C$ is a uniformly bounded constant. By defining $\eps_\beta =
(1+C\eps_\alpha)(1+\eps)-1$, we have
\[
\lrsa(\beta) = \beta (1+\eps_\beta), \quad |\eps_\beta|\lesssim \eps.
\]

For $\gamma = i\omega(x_1-y_1)\beta$, we define $\lrsa(\gamma) =
i\omega(x_1-y_1) \lrsa(\beta)$ and $\eps_\gamma =
\eps_\beta$. Clearly,
\[
\lrsa(\gamma) = \gamma (1+\eps_\gamma), \quad |\eps_\gamma|\lesssim \eps.
\]
An easy but essential calculation also shows that
\[
\gamma = i\omega(x_1-y_1) \left(\sqrt{1+\frac{|x_2-y_2|^2}{|x_1-y_1|^2}}-1\right)
\]
is of order $O(1)$ for $x\in R_1$ and $y\in R_2$. Then it follows that
$\lrsa(\gamma)$ is also of order $O(1)$.

Finally, consider $\eta = \exp(i\gamma)$. As both $\gamma$ and
$\lrsa(\gamma)$ are of order $O(1)$, applying truncated Taylor
expansion to $\exp(i\gamma)$ near $\gamma=0$ with relative error
$\eps$ gives a polynomial $q(\gamma)$ with an $(\omega,d)$-independent
number of terms. We then define $\lrsa(\eta) = q(\lrsa(\gamma))$ and
have
\begin{align*}
  \lrsa(\eta) &= q(\lrsa(\gamma)) = \exp(i\lrsa(\gamma)) (1+\eps) \\
  &= \exp(i\gamma) \exp(i(\lrsa(\gamma)-\gamma)) (1+\eps) \\
  &= \exp(i\gamma) \exp(i\gamma \eps_\gamma) (1+\eps) \\
  &= \exp(i\gamma) (1+C\gamma \eps_\gamma) (1+\eps).
\end{align*}
where $C$ is uniformly bounded from the fact that $\gamma=O(1)$ and
$\lrsa(\gamma)=O(1)$. By defining $\eps_\eta = (1+C\gamma
\eps_\gamma) (1+\eps)-1$, one has
\[
\lrsa(\eta) = \eta (1+\eps_\eta), \quad |\eps_\eta|\lesssim \eps.
\]
This shows that $\exp(i\omega|x-y|)$ has a low-rank separated
approximation with an $(\omega,d)$-independent number of terms and
relative error $O(\eps)$. Combined with the approximation for
$G_0(|x-y|)$, we conclude that $G(x,y)$ has a low-rank separated
approximation with an $(\omega,d)$-independent number of terms and
relative error $O(\eps)$.

At this point, it is clear that in order to get a relative error
strictly bounded by $\eps$, we only need to repeat the above argument
by pre-dividing $\eps$ by a uniform constant factor.
\end{proof}

\subsection{Setup algorithm}


We call that the length of $\p\Omega$ is $4^q\lambda$. The setup
algorithm first generates a complete binary tree structure of
$\p\Omega$ by bisecting it into segments of equal length recursively.
The process is stopped when a segment is of length $m_\ell \lambda$,
where $m_\ell$ is a constant typically set to 2 or 4. A segment at the
final level is called {\em a leaf} and it contains $m_\ell p$
discretization points. A segment $T$ is called {\em almost-planar} if
its length is bounded by by $2^q \lambda/\sqrt{\kappa_T}$, where
$\kappa_T$ is the maximum of the absolute value of the curvature
within the segment $T$.

\begin{figure}[h!]
  \begin{center}
    \includegraphics[height=1in]{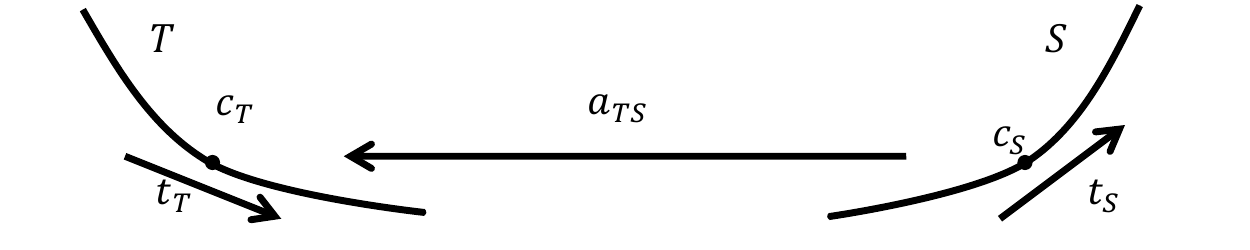}    
  \end{center}
  \caption{The geometric setting for constructing low-rank
    approximation between two segments $T$ and $S$.}
  \label{fig:dir}
\end{figure}

The key component of the setup algorithm is how to compress the
interaction between two segments. Consider two segments $T$ and $S$
and let
\begin{itemize}
\item $c_T$ and $c_S$ be the centers of $T$ and $S$, 
\item $t_T$ and $t_S$ be the tangent directions at $c_T$ and $c_S$, and 
\item $a_{TS}=(c_T-c_S)/|c_T-c_S|$ be the unit vector pointing from $c_S$
to $c_T$ (see Figure \ref{fig:dir}).
\end{itemize}

\newcommand{\dist}{\text{dist}}
\newcommand{\diam}{\text{diam}}

The pair $(T,S)$ is called {\em parabolically separated} if $T$ and
$S$ can be enclosed in two rectangles $R_1$ and $R_2$ that satisfy the
condition of Theorem \ref{thm:main}. More precisely, we define
\begin{itemize}
\item $\ell_{TS}$ to be the line passing through the centers $c_T$
  and $c_S$,
\item $w_T\lambda$ and $w_S\lambda$ to be the lengths of the
  projections of $T$ and $S$ onto $\ell_{TS}$,
\item $d_{TS}\lambda$ to be the distance between the projections, and
\item $h_T\lambda$ and $h_S\lambda$ to be the lengths of the
  projections of $T$ and $S$, respectively, onto the direction
  $a_{TS}^\perp$.
\end{itemize}
With these definitions, $(T,S)$ is parabolically separated if
$d_{TS}>1$, $d_{TS}>\max(w_T,w_S)$, and
$d_{TS}>2\max(h_T,h_S)^2$. Theorem \ref{thm:main} states that, for a
parabolically separated pair $(T,S)$, the interaction between $T$ and
$S$ is numerically low-rank. As we shall see, the setup algorithm
constructs a low-rank approximation for such a pair when at least one
of $T$ and $S$ is also almost-planar.

Let $\G$ be an empty set initially. The algorithm first traverses the
tree in a breadth first search. Whenever it reaches an almost-planar
segment, the algorithm inserts it into $\G$ and skip the whole subtree
starting from it. When it reaches a leaf, the algorithm also inserts
it into $\G$. At the end of this step, $\G$ consists of almost-planar
segments and non-planar leaves, and they form a disjoint union of the
boundary $\p\Omega$.

Next, the algorithm starts with a queue $\Q$ of segment pairs,
initialized to be $\{(T,S), T,S\in \G\}$, as well as a set $\F$,
initialized to be empty. When the queue $\Q$ is not empty, the top
element is popped out and denoted by $(T,S)$. The algorithm treats the
pair $(T,S)$ depending on the types of $T$ and $S$.

\begin{itemize}
\item If both $T$ and $S$ are non-planar leaves, we mark $(T,S)$
  as a {\em dense pair} and insert it into $\F$.
\item Suppose that $T$ is a non-planar leaf but $S$ is
  almost-planar. If $(T,S)$ is parabolically separated, then we mark
  $(T,S)$ as a {\em low-rank pair} and insert it into $\F$. If it is not,
  one of the following two options is taken.
  \begin{itemize}
  \item If $S$ is not at the leaf level, we partition $S$ evenly into
    two segments $S_1$ and $S_2$ and insert $(T,S_1)$ and $(T,S_2)$
    into $\Q$.
  \item If $S$ is at the leaf level, we mark $(T,S)$ as a {\em dense
    pair} and insert it into $\F$.
  \end{itemize}
\item Suppose that $T$ is almost-planar but $S$ is a non-planar
  leaf. If $(T,S)$ is parabolically separated, then we mark $(T,S)$ as
  a {\em low-rank pair} and insert it into $\F$. If it is not, one of the
  following two options is taken.
  \begin{itemize}
    \item If $T$ is not at the leaf level, we partition $T$ evenly
      into two segments $T_1$ and $T_2$ and insert $(T_1,S)$ and
      $(T_2,S)$ into $\Q$. 
    \item If $T$ is at the leaf level, we mark $(T,S)$ as a {\em dense
      pair} and insert it into $\F$.
  \end{itemize}
\item Suppose that both $T$ and $S$ are almost-planar. If $(T,S)$ is
  parabolically separated, then we mark $(T,S)$ as a {\em low-rank
    pair} and insert it into $\F$. If it is not, one of the following
  options is taken.
  \begin{itemize}
  \item If $T$ is longer than $S$, partition $T$ evenly into two
    segments $T_1$ and $T_2$ and insert $(T_1,S)$ and $(T_2,S)$ into
    $\Q$. 
  \item If $S$ is longer than $T$, partition $S$ evenly into two
    segments $S_1$ and $S_2$ and insert $(T,S_1)$ and $(T,S_2)$ into
    $\Q$.  
  \item If $T$ and $S$ are of the same length but are not at the leaf
    level, then partition $T$ evenly into $T_1$ and $T_2$ and $S$
    evenly into $S_1$ and $S_2$ and insert $(T_1,S_1)$, $(T_2,S_1)$,
    $(T_1,S_2)$, and $(T_2,S_2)$ into $\Q$. 
  \item Finally if $T$ and $S$ are of the same length and at the leaf
    level, then mark $(T,S)$ as a {\em dense pair} and insert it into
    $\F$.
  \end{itemize}
\end{itemize}
This process is repeated as long as the queue $\Q$ is not empty. When
$Q$ is empty, the set $\F$ provides a hierarchical compressed
approximation of \eqref{eq:uxdis}.

\subsubsection{Construction of directional low-rank approximation}
What is missing at this point is the procedure of constructing a
actual low-rank separated approximation. In the setup algorithm, this
procedure is used when $(T,S)$ is parabolically separated and at least
one of $T$ and $S$ is almost-planar.

Let us first consider a low-rank pair $(T,S)\in\F$ where both $T$ and
$S$ are almost-planar. It is instructive to consider the function
$\exp(i\omega|x-y|)$ instead of $G(x,y)$ as this is where the
oscillations come from. In what follows, we shall use $\sim$ to denote an
approximation up to a non-oscillatory multiplicative term.

Let $2^{\ell_T}\lambda$ and $2^{\ell_S}\lambda$ be the lengths of $T$
and $S$, respectively. For a point $x\in T$ and $y\in S$,
\begin{align}
  &\exp(i\omega|x-y|) \sim \exp(i\omega  a_{TS}\cdot(x-y)) \nonumber \\
  =&\exp(i\omega a_{TS}\cdot((x-c_T)+(c_T-c_S)+(c_S-y)) ) \nonumber \\
  =&\exp(i\omega a_{TS}\cdot(x-c_T)) \cdot \exp(i\omega a_{TS}\cdot (c_T-c_S)) 
  \cdot \exp(-i\omega a_{TS}\cdot(y-c_S)), \label{eq:expsep}
\end{align}
where the approximation in the first line is based on the fact that
$(T,S)$ is parabolically separated. 

To further approximate the first term in \eqref{eq:expsep}, we perform
a Taylor expansion for $\rho^{-1}(t)$ near $t=\rho(c_T)$ and evaluate
it at $\rho(x)$:
\begin{align*}
  & \left| \rho^{-1}(\rho(x)) - \left( \rho^{-1}(\rho(c_T)) + t_T (\rho(x)-\rho(c_T)) \right) \right|\\
  & \lesssim \frac{1}{2} |\rho(x)-\rho(c_T)|^2 \kappa_T \le \frac{1}{2} (2^q\lambda/\sqrt{\kappa_T})^2 \kappa_T
  = O(4^q \lambda^2)  = O(\lambda),
\end{align*}
where $\kappa_T$ is the maximum of the absolute value of the curvature
in $T$ and here we also use the fact that $T$ is almost-planar. This
is equivalent to
\[
(x-c_T) = (\rho(x)-\rho(c_T))\cdot t_T + O(\lambda).
\]
Multiplying it with $i \omega a_{TS}$ and taking exponential gives the approximation
\begin{equation}
  \exp( i\omega a_{TS}\cdot(x-c_T) ) \sim \exp(i \omega a_{TS}\cdot t_T (\rho(x)-\rho(c_T))).
\label{eq:tflat}
\end{equation}

For the last term of \eqref{eq:expsep}, the same argument works for
$(y-c_S)$. Since $S$ is almost-planar, we get
\[
(y-c_S) = (\rho(y)-\rho(c_S))\cdot t_S + O(\lambda),
\]
and 
\begin{equation}
  \exp(-i\omega a_{TS}\cdot (y-c_S)) \sim \exp(-i \omega a_{TS}\cdot t_S (\rho(y)-\rho(c_S))).
\label{eq:sflat}
\end{equation}

Noticing that $\omega a_{TS}\cdot t_T \in [-\omega,\omega]$, we
partition the interval $[-\omega,\omega]$ uniformly into
$2^{\ell_T+1}m_f$ intervals with a set $K_T$ of $2^{\ell_T+1}m_f+1$
gridpoints. Here $m_f$ is a parameter that is typically set to $2$ or
$4$. A key step of the algorithm is to approximate $\omega a_{TS}\cdot
t_T$ with a nearby gridpoint. More precisely, we define $[k]_T$ to the
value of rounding $k$ to the nearest gridpoint in $K_T$. Since
\[
\left|(\omega a_{TS}\cdot t_T - [\omega a_{TS}\cdot t_T]_T)(\rho(x)-\rho(c_T)) \right| \le
\frac{2\omega}{2^{\ell_T+1}m_f} \cdot \frac{1}{2}\cdot \frac{2^{\ell_T}\lambda}{2} = \frac{2\pi}{4 m_f} = O(1),
\]
replacing $\omega a_{TS}\cdot t_T$ with $[\omega a_{TS}\cdot t_T]_T$
in \eqref{eq:tflat} only introduces a non-oscillatory multiplicative
term. Therefore,
\begin{equation}
  \exp( i\omega a_{TS}\cdot(x-c_T) ) \sim \exp(i [\omega a_{TS}\cdot t_T]_T (\rho(x)-\rho(c_T))).
  \label{eq:tflatnew}
\end{equation}

Similarly for $-\omega a_{TS}\cdot t_S$, we partition the interval
$[-\omega,\omega]$ uniformly into $2^{\ell_S+1}m_f$ intervals with a
set $K_S$ of $2^{\ell_S+1}m_f+1$ gridpoints and define $[k]_S$ to the
value of rounding $k$ to the nearest gridpoint in $K_S$. Again since
\[
\left|(-\omega a_{TS}\cdot t_S - [-\omega a_{TS}\cdot t_S]_S)(\rho(y)-\rho(c_S)) \right| \le
\frac{2\omega}{2^{\ell_S+1}m_f} \cdot \frac{1}{2}\cdot \frac{2^{\ell_S}\lambda}{2} = \frac{2\pi}{4 m_f} = O(1),
\]
replacing $-\omega a_{TS}\cdot t_S$ with $[-\omega a_{TS}\cdot t_S]_S$
in \eqref{eq:sflat} introduces an extra non-oscillatory multiplicative
term. Thus
\begin{equation}
  \exp(-i\omega a_{TS}\cdot (y-c_S)) \sim \exp(i [-\omega a_{TS}\cdot t_S]_S (\rho(y)-\rho(c_S))).
  \label{eq:sflatnew}
\end{equation}
The reason for introducing the rounding is that: 
\begin{itemize}
\item $\exp(i k (\rho(x)-\rho(c_T))$ for $x\in T$ and $k\in K_T$ is a
  partial Fourier matrix since both $\rho(x)$ and $k$ are on a uniform
  grid.
\end{itemize}
The same also holds for $\exp(i k (\rho(y)-\rho(c_S))$ for $y\in S$
and $k\in K_S$.  As we shall see, this has an enormous impact on the
efficiency of applying the low-rank approximations.

By introducing 
\begin{align*}
  &k^T_{TS} = [\omega a_{TS} t_T]_T,\\
  &k^S_{TS} = [-\omega a_{TS} t_S]_S,\\
  &F_T(x,k) = \exp(i k (\rho(x)-\rho(c_T))), \quad x\in T, k\in K_T\\
  &F_S(y,k) = \exp(i k (\rho(y)-\rho(c_S))), \quad y\in S, k\in K_S
\end{align*}
and applying \eqref{eq:tflatnew} and \eqref{eq:sflatnew} to
\eqref{eq:expsep}, we have
\[
\exp(i\omega|x-y|) \sim F_T(x,k^T_{TS}) \cdot
\exp(i\omega(c_T-c_S)\cdot a_{TS}) \cdot F_S(y,k^S_{TS}).
\]

So far, we have been considering the kernel $\exp(i\omega|x-y|)$. Let
us now replace $\exp(i\omega|x-y|)$ back with the kernel
$G(x,y)$. Since the difference $G_0(|x-y|)$ is a non-oscillatory term,
the above discussion applies without any change and we have the
following representation for $G(x,y)$
\begin{equation}
  G(x,y) = F_T(x,k^T_{TS}) \cdot G^{mm}_{TS}(x,y) \cdot F_S(y,k^S_{TS}),\quad x\in T, y\in S,
  \label{eq:Gxy}
\end{equation}
where the term $G^{mm}_{TS}(x,y)$ defined via this equation is
non-oscillatory in $x\in T$ and $y\in S$. The superscript $mm$ of
$G^{mm}_{TS}$ indicates that the complex exponential modulation is
done both for $T$ and for $S$. In a more compact operator form, we can
write \eqref{eq:Gxy} as
\begin{equation}
  G(T,S) = \diag(F_T(:,k^T_{TS})) \cdot G^{mm}_{TS}(T,S) \cdot \diag(F_S(:,k^S_{TS})).
  \label{eq:Gxyop}
\end{equation}

Since $G^{mm}_{TS}(x,y)$ is non-oscillatory, it can be approximated
with Chebyshev interpolation in the parametric domain of the
boundary. More specifically, we define
\begin{itemize}
\item $m_c$ to be the size of the Chebyshev
  grid used,
\item $R_T$ be the image (under $\rho^{-1}$) of the Chebyshev grid in
  $T$ and $R_S$ be the image of the Chebyshev grid in $S$,
\item $I_T$ and $I_S$ to be the interpolation operator for $T$ and $S$
  associated with the Chebyshev grids $R_T$ and $R_S$, respectively.
\end{itemize}
In the matrix notation, $I_T$ is a matrix with entries given by
$I_T(x,j)$ for $x\in T$ and $j\in R_T$ and $I_S$ is a matrix with
entries given by $I_S(y,j)$ for $y\in S$ and $j\in R_S$.  In an
operator form, this approximation reads
\[
G^{mm}_{TS}(T,S) \approx I_T \cdot G^{mm}_{TS}(R_T,R_S) \cdot I_S^t.
\]
Putting it together with \eqref{eq:Gxyop} gives
\[
G(T,S) \approx \diag(F_T(:,k^T_{TS})) \cdot I_T \cdot
G^{mm}_{TS}(R_T,R_S) \cdot I_S^t \cdot \diag(F_S(:,k^S_{TS})).
\]
Notice that for a pair $(T,S)$, representing this low-rank
approximation only requires storing $k^T_{TS}$, $k^S_{TS}$, and
$G^{mm}_{TS}(R_T,R_S)$.

The above discussion is for the case in which both $T$ and $S$ are
almost-planar. However, the setup algorithm also needs to construct
low-rank approximations for low-rank pairs in which only one of $T$
and $S$ is almost-planar. For those cases, the oscillatory Chebyshev
interpolation is only used for the almost-planar segment while a dense
evaluation is used at the other non-planar segment. More precisely, if
$S$ is almost-planar, we form the representation
\[
G(T,S) = G^{\cdot m}_{TS}(T,S) \cdot \diag(F_S(:,k^S_{TS}))
\]
and the approximation
\[
G(T,S) \approx G^{\cdot m}_{TS}(T,R_S) \cdot I_S^t \cdot  \diag(F_S(:,k^S_{TS})),
\]
where the superscript $\cdot m$ means that the exponential modulation
is done only for $S$. Similarly if $T$ is almost-planar, we construct
the representation
\[
G(T,S) = \diag(F_T(:,k^T_{TS})) \cdot G^{m\cdot}_{TS}(T,S)
\]
and the approximation
\[
G(T,S) \approx \diag(F_T(:,k^T_{TS})) \cdot I_T \cdot G^{m\cdot}_{TS}(R_T,S),
\]
where the superscript $m\cdot$ means that the exponential modulation
is done only for $T$.

\subsubsection{Complexity of the setup algorithm}

Since $\Omega$ has a $C^2$ boundary, the curvature is uniformly
bounded. Let us first estimate the number of segments in $\G$. From
the setup algorithm, it is clear that the length of a segment $T$ in
$\G$ is bounded by $\min(2^q\lambda, 2^q\lambda/\sqrt{\kappa_T}) =
\Theta(2^q\lambda) = \Theta(\sqrt{\omega}\lambda)$. As the total
length of the boundary is $\Theta(1)$, the number of segments in $\G$
is bounded by $\Theta(1/(\sqrt{\omega}\lambda)) =
\Theta(\sqrt{\omega})$.

Next, let us consider the number of pairs in $\F$. Since that the
boundary is $C^2$, for almost all pairs $(T,S)$ in $\F$ the segments
$T$ and $S$ have the same length. Consider first the segments of
length $2^q\lambda$. As there are $2^q$ segments of this length, the
total number of pairs is bounded by $2^{2q} = \Theta(\omega)$. For
length equal to smaller values, a segment $T$ only appears together
with a constant number of segments $S$ in its neighborhood. Therefore,
the total number of pairs in $\F$ is again $\Theta(\omega)$.

The complexity of the setup algorithm contains three parts:
\begin{itemize}
\item The generation of the set $\F$. It takes $O(\omega\log\omega)$
  steps.
\item The evaluation of the dense pairs. Clearly there are at most
  $\Theta(\omega)$ dense pairs. Since each leaf segment contains at
  most $O(1)$ points, the dense matrix $G(T,S)$ has $O(1)$
  entries. Therefore, the evaluation cost of all dense pairs is
  $O(\omega)$.
\item The evaluation of the low-rank pairs. Again, the number of
  low-rank pairs is at most $\Theta(\omega)$. In all cases of the
  low-rank pairs, the matrices $G^{mm}_{T,S}$, $G^{\cdot m}_{T,S}$,
  and $G^{m\cdot }_{T,S}$ have $O(1)$ entries.  Therefore, the
  evaluation cost of all low-rank pairs is also $O(\omega)$.
\end{itemize}
Summing these contributions together gives a $\Theta(\omega\log\omega)
= \Theta(n\log n)$ complexity for the setup algorithm.

\subsection{Evaluation algorithm}

Now we are ready to describe how to evaluate the sums
\[
u(x) = \sum_{y\in P} G(x,y) f(y), \quad x\in P
\]
efficiently. In the following discussion, we slightly abuse notation
by using $u$ to denote the vector $(u(x))_{x\in P}$ and $f$ to denote
the vector $(f(y))_{y\in P}$. For a segment $T$, $u(T)$ stands for the
subvector obtained by restricting $u$ to the points in $T$. Similarly,
$f(S)$ is the restriction of $f$ to the points in $S$.

Initially, we set $u=0$. The algorithm visits all pairs $(T,S)\in
\F$. If $(T,S)$ is dense, then
\[
u(T) \Leftarrow u(T) + G(T,S) f(S),
\]
where $\Leftarrow$ stands for assignment. 

If $(T,S)$ is low-rank and both $T$ and $S$ are almost-planar,
\begin{equation}
  u(T) \Leftarrow u(T) + \diag(F_T(:,k^T_{TS})) \cdot I_T \cdot
  G^{mm}_{TS}(R_T,R_S) \cdot I_S^t \cdot \diag(F_S(:,k^S_{TS}))
  \cdot f(S).
  \label{eq:Gmm}
\end{equation}
If $(T,S)$ is low-rank and only $S$ is almost-planar,
\begin{equation}
  u(T) \Leftarrow u(T) + G^{\cdot m}_{TS}(T,R_S) \cdot I_S^t \cdot
  \diag(F_S(:,k^S_{TS})) \cdot f(S).
  \label{eq:Gdm}
\end{equation}
If $(T,S)$ is low-rank and only $T$ is almost-planar,
\begin{equation}
  u(T) \Leftarrow u(T) + \diag(F_T(:,k^T_{TS})) \cdot I_T \cdot
  G^{m\cdot}_{TS}(R_T,S) \cdot f(S).
  \label{eq:Gmd}
\end{equation}

Direct evaluation of these three formulas results in a
$\Theta(\omega^{3/2})$ complexity. In order speed up this calculation,
the work is split into multiple steps. Let us consider \eqref{eq:Gmm},
i.e., the situation where both $T$ and $S$ are almost-planar.

First, at segment $S$, one needs to compute
\[
I_S^t \cdot \diag(F_S(:,k^S_{TS})) \cdot f(S).
\]
It turns out to be extremely useful to consider all $k \in K_S$
altogether instead of one by one:
\[
I_S^t \cdot \diag(F_S(:,k)) \cdot f(S).
\]
This is equivalent to evaluating for each $j \in R_S$
\begin{equation}
  \langle I_S(:,j), F_S(:,k) \odot f(S) \rangle =
  \langle F_S(:,k), I_S(:,j) \odot f(S) \rangle := \hat{f}_S(j,k),
  \label{eq:fhat}
\end{equation}
where the symbol $\odot$ stands for the entrywise product of two
vectors and $\langle\cdot,\cdot\rangle$ is the bilinear product, i.e.,
$\langle v, w \rangle = \sum_{i} v_i w_i$. An essential observation
from the last equation is that:
\begin{itemize}
\item Computing $\hat{f}_S(j,:)$ (i.e., all $k\in K_S$ with a fixed
  $j\in R_S$) simply requires entrywise product of $f(S)$ with
  $I_S(:,k)$, followed by an FFT.
\end{itemize}
Looping through each $j\in R_S$ provides all entries in the
$|R_S|\times|K_S|$ matrix $\hat{f}_S(:,:)$.

Second, instead of considering one pair $(T,S)$ at a time, we fix $T$ and
consider all pairs for which the first component is $T$ and both $T$
and $S$ are almost-planar. More precisely,
\begin{align*}
u(T) &\Leftarrow u(T) + \sum_S \diag(F_T(:,k^T_{TS})) \cdot I_T \cdot G^{mm}_{TS}(R_T,R_S) \cdot 
\hat{f}_S(:,k^S_{TS})\\
&\Leftarrow u(T) + \sum_{k\in K_T} \sum_{S: k^T_{TS}=k} \diag(F_T(:,k)) \cdot I_T \cdot G^{mm}_{TS}(R_T,R_S) \cdot 
\hat{f}_S(:,k^S_{TS})\\
&\Leftarrow u(T) + \sum_{k\in K_T} \diag(F_T(:,k)) \cdot I_T \cdot \left(\sum_{S: k^T_{TS}=k} G^{mm}_{TS}(R_T,R_S) \cdot 
\hat{f}_S(:,k^S_{TS})\right).
\end{align*}
Motivated by the last equation, we introduce an $|R_T|\times|K_T|$
matrix $\hat{u}_T(:,:)$ defined via
\[
\hat{u}_T(:,k) = \sum_{S: k^T_{TS}=k} G^{mm}_{TS}(R_T,R_S) \hat{f}_S(:,k^S_{TS}).
\]

Third, we can now write
\[
u(T) \Leftarrow u(T) + \sum_{k \in K_T} \diag(F_T(:,k)) \cdot I_T \cdot \hat{u}_T(:,k).
\]
To perform this computation efficiently, notice that the last term is
\begin{equation}
  \sum_{k\in K_T} \sum_{j\in R_T} \diag(F_T(:,k)) I_T(:,j) \cdot \hat{u}_T(j,k)
  = \sum_{j\in R_T} I_T(:,j) \odot \left(\sum_{k\in K_T} F_T(:,k) \hat{u}_T(j,k) \right).
  \label{eq:uhat}
\end{equation}

An essential observation is that
\begin{itemize}
\item for a fixed $j\in R_T$, the sum over
  $k\in K_T$ can be carried out by an FFT. 
\end{itemize}
Therefore, for each $j\in R_T$, the work is reduced to an FFT followed
by an entrywise product.

The above discussion addresses the case where both $T$ and $S$ are
almost-planar. For the cases where only one of them is almost-planar,
i.e., \eqref{eq:Gdm} and \eqref{eq:Gmd}, the procedure is similar
except that the oscillatory Chebyshev interpolation and local FFTs are
done only on one side.

To summarize, the evaluation algorithm takes the following steps:
\begin{enumerate}
\item Set $u=0$ and $\hat{u}_T(:,:) = 0$.
\item For each segment $S$ and for each $j\in R_S$, form
  $\hat{f}_S(j,:)$ by an entrywise product followed by an FFT, as shown
  in \eqref{eq:fhat}.
\item For each pair $(T,S)$ in $\F$,
  \begin{enumerate}
    \item If $(T,S)$ is dense,
      \[
      u(T) \Leftarrow u(T) + G(T,S) f(S).
      \]
    \item If $(T,S)$ is low-rank and both are almost-planar,
      \[
      \hat{u}_T(:,k^T_{TS}) \Leftarrow \hat{u}_T(:,k^T_{TS}) + G^{mm}_{TS}(R_T,R_S) \hat{f}_S(:,k^S_{TS}).
      \]
    \item If $(T,S)$ is low-rank and only $S$ is almost-planar,
      \[
      u(T) \Leftarrow u(T) + G^{\cdot m}_{TS}(T,R_S) \hat{f}_S(:,k^S_{TS}).
      \]
    \item If $(T,S)$ is low-rank and only $T$ is almost-planar,
      \[
      \hat{u}_T(:,k^T_{TS}) \Leftarrow \hat{u}_T(:,k^T_{TS}) + G^{m\cdot}_{TS}(R_T,S) f(S).
      \]
  \end{enumerate}
\item For each segment $T$ and for each $j\in R_T$, apply an FFT and
  then an entrywise product to $\hat{u}_T(j,:)$ and add the result to
  $u(T)$, as shown in \eqref{eq:uhat}.
\end{enumerate}

\subsubsection{Complexity of the evaluation algorithm}

The complexity of the four steps of the evaluation algorithm is
estimated as follows.
\begin{enumerate}
\item The cost of the first step is clearly $O(\omega)$.
\item The FFT for a segment of length $t\lambda$ and $tp$ points take
  $\Theta(t\log t)$ steps as $p$ is constant. Since there are
  $O(4^q/t)$ segments of this length, the overall cost is $O(4^q \log
  t) = O(\omega \log t)$. Now summing over all possible values of $t$
  results a complexity of $O(\omega \log^2 \omega)$ for the second
  step of the algorithm.
\item Now consider the third step. The matrices $G(T,S)$,
  $G^{mm}_{TS}(R_T,R_S)$, $G^{\cdot m}_{TS}(T,R_S)$, and
  $G^{m\cdot}_{TS}(R_T,S)$ that appear in all four scenarios all have
  size $O(1)$.  Taking into consideration that there are at most
  $O(\omega)$ pairs in $\F$, the cost of this step is $O(\omega)$.
\item The last step is very similar to the second one. Its cost is
  $O(\omega\log^2\omega)$ as well.
\end{enumerate}
Putting these together shows that the complexity of the evaluation algorithm is
$O(\omega \log^2\omega) = O(n\log^2 n)$.

\section{Numerical Results}

\begin{figure}[h!]
  \begin{center}
    \includegraphics[height=1.5in]{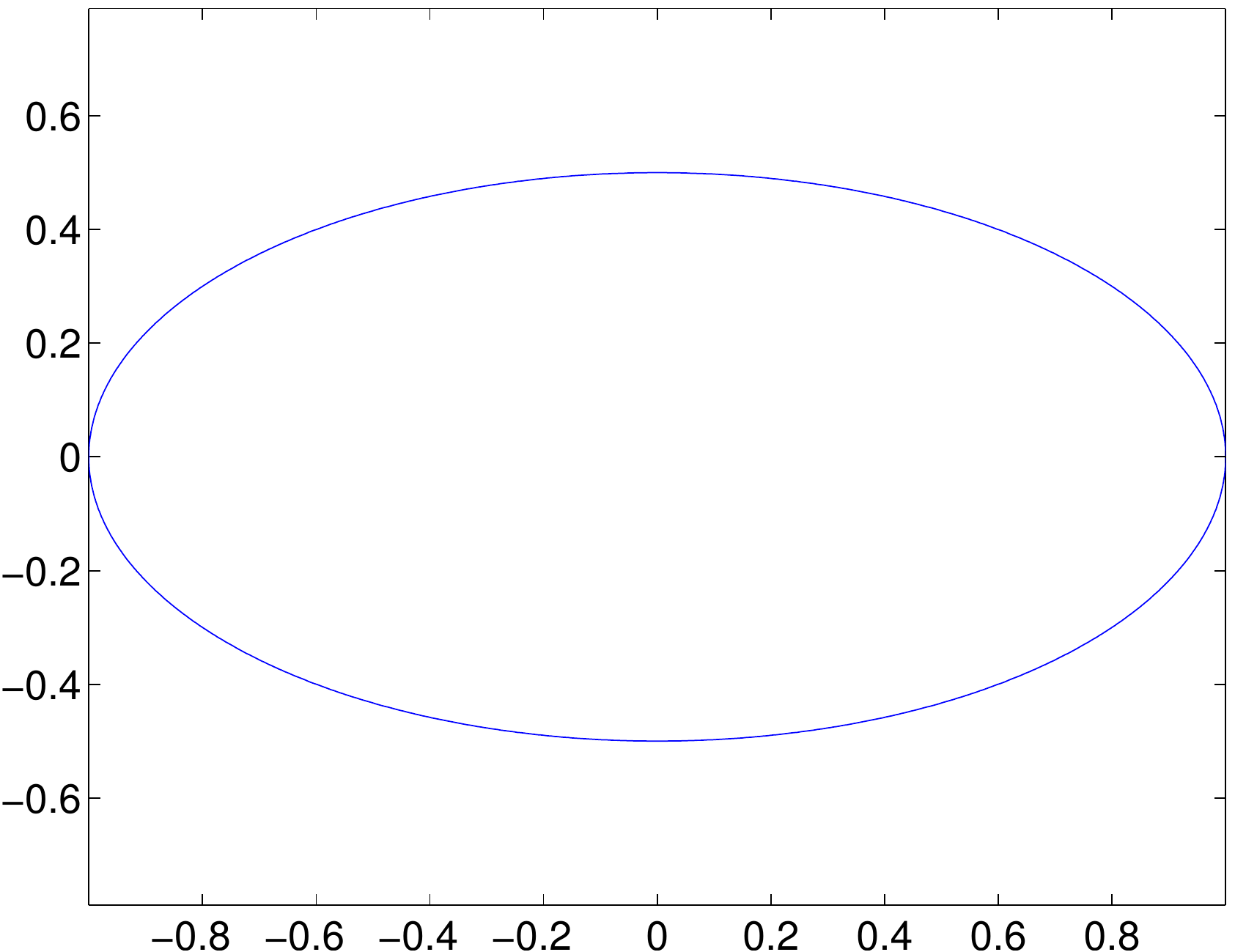}    \hspace{0.25in}
    \includegraphics[height=1.5in]{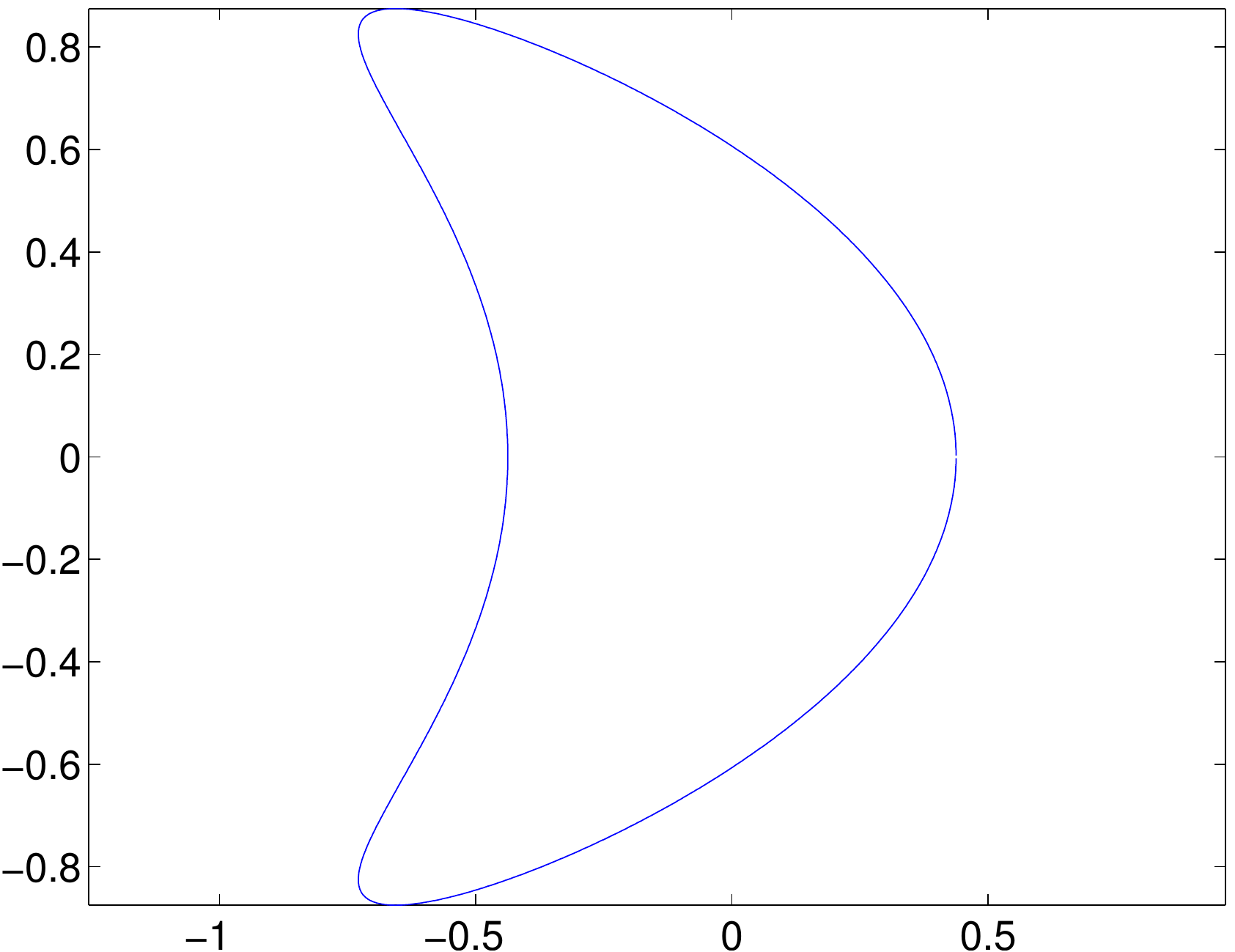}    
  \end{center}
  \caption{The two scatterers used in the numerical tests. (a) an
    ellipse. (b) a bean-shaped object.}
  \label{fig:doms}
\end{figure}

The proposed algorithms are implemented in Matlab and the numerical
results in this section are obtained on a desktop computer with a
3.60GHz CPU. The numerical experiments are performed for two domains:
(a) an ellipse and (b) a bean-shaped object (shown in Figure
\ref{fig:doms}).

The constants in the algorithm are set as follows. 
\begin{itemize}
\item The number of points $p$ per wavelength is $8$.
\item The length of the leaf segment is $m_\ell\lambda$ with $m_\ell=4$.
\item The constant $m_f$ in defining the frequency grid is equal to
  $2$.
\item The size $m_c$ of the Chebyshev grid is equal to $6$, $8$, $10$,
  or $12$. This number controls the accuracy of the evaluation
  algorithm.
\end{itemize}

We denote the exact values and the numerical approximations by
$u_e(x)$ and $u_a(x)$, respectively. To check the error, the exact
solution is evaluated at a set $S$ of 100 points and the relative error
is estimated by
\begin{equation}
  \left( 
  \frac{\sum_{x\in S} |u_e(x)-u_a(x)|^2}{\sum_{x\in S} |u_e(x)|^2}
  \right)^{1/2}.
  \label{eq:err}
\end{equation}

The experiments are performed for the following four operators
\begin{align*}
  u(x) &= (Sf)(x) := \int_{\p\Omega} G(x,y) f(y) dy,\\
  u(x) &= (Df)(x) := \int_{\p\Omega} \frac{\p G(x,y)}{\p n(y)} f(y) dy,\\
  u(x) &= (D'f)(x) := \int_{\p\Omega} \frac{\p G(x,y)}{\p n(x)} f(y) dy,\\
  u(x) &= (Nf)(x) := \int_{\p\Omega} \frac{\p^2 G(x,y)}{\p n(x)\p n(y)} f(y) dy,
\end{align*}
which are typically referred as single layer ($S$), double layer
($D$), normal derivative of single layer ($D'$), and normal derivative
of double layer ($N$) \cite{colton-2013,nedelec-2001}. These operators
are the basic building blocks of the boundary integral methods for
solving the Dirichlet and Neumann problems of acoustic wave
scattering.

The numerical results for the single layer operator $S$ on the two
domains are given in Tables \ref{tbl:S1} and \ref{tbl:S2}. In these
tables,
\begin{itemize}
\item $m_c$ is the size of the Chebyshev grid that controls the
  accuracy of the method,
\item $\omega$ is the frequency of the wave field,
\item $n$ is the number of boundary discretization points,
\item $T_s$ is the running time for the setup algorithm,
\item $T_a$ is the running times of the application algorithm, and
\item $e$ is the relative error estimated via \eqref{eq:err}.
\end{itemize}
The numerical results for the double layer operator $D$ on the two
domains are given in Tables \ref{tbl:D1} and \ref{tbl:D2}.  The
numerical results for the normal derivative of the single layer
operator $D'$ on the two domains are given in Tables \ref{tbl:X1} and
\ref{tbl:X2}. Finally, the numerical results for the normal derivative
of the double layer operator $N$ on the two domains are given in
Tables \ref{tbl:N1} and \ref{tbl:N2}.


\begin{table}[ht!]
  \begin{center}
    \begin{tabular}{|c|cc|ccc|}
      \hline
      $m_c$ & $\omega$ & $n$ & $T_s$ & $T_a$ & $e$\\
      \hline
      6 & 5.3e+03 & 3.3e+04 & 1.9e+01 & 1.3e+00 & 2.4e-04\\
      6 & 2.1e+04 & 1.3e+05 & 8.0e+01 & 5.5e+00 & 3.5e-04\\
      6 & 8.5e+04 & 5.2e+05 & 3.3e+02 & 2.2e+01 & 2.9e-04\\
      \hline
      8 & 5.3e+03 & 3.3e+04 & 2.0e+01 & 1.9e+00 & 1.0e-05\\
      8 & 2.1e+04 & 1.3e+05 & 8.2e+01 & 5.6e+00 & 8.8e-06\\
      8 & 8.5e+04 & 5.2e+05 & 3.3e+02 & 2.4e+01 & 1.3e-05\\
      \hline
      10 & 5.3e+03 & 3.3e+04 & 2.1e+01 & 2.0e+00 & 9.9e-07\\
      10 & 2.1e+04 & 1.3e+05 & 8.4e+01 & 6.1e+00 & 7.8e-07\\
      10 & 8.5e+04 & 5.2e+05 & 3.4e+02 & 2.6e+01 & 5.3e-07\\
      \hline
      12 & 5.3e+03 & 3.3e+04 & 2.2e+01 & 2.2e+00 & 2.0e-08\\
      12 & 2.1e+04 & 1.3e+05 & 8.8e+01 & 6.6e+00 & 2.2e-08\\
      12 & 8.5e+04 & 5.2e+05 & 3.6e+02 & 2.9e+01 & 2.7e-08\\
      \hline
    \end{tabular}
  \end{center}
  \caption{Numerical results of operator $S$ for the ellipse.}
  \label{tbl:S1}
\end{table}

\begin{table}[ht!]
  \begin{center}
    \begin{tabular}{|c|cc|ccc|}
      \hline
      $m_c$ & $\omega$ & $n$ & $T_s$ & $T_a$ & $e$\\
      \hline
      6 & 5.2e+03 & 3.3e+04 & 2.6e+01 & 2.0e+00 & 7.6e-04\\
      6 & 2.1e+04 & 1.3e+05 & 1.0e+02 & 6.6e+00 & 9.1e-04\\
      6 & 8.3e+04 & 5.2e+05 & 4.2e+02 & 3.0e+01 & 9.4e-04\\
      \hline
      8 & 5.2e+03 & 3.3e+04 & 2.7e+01 & 2.4e+00 & 5.6e-05\\
      8 & 2.1e+04 & 1.3e+05 & 1.1e+02 & 7.0e+00 & 5.4e-05\\
      8 & 8.3e+04 & 5.2e+05 & 4.3e+02 & 3.0e+01 & 7.5e-05\\
      \hline
      10 & 5.2e+03 & 3.3e+04 & 2.8e+01 & 2.5e+00 & 2.2e-06\\
      10 & 2.1e+04 & 1.3e+05 & 1.1e+02 & 7.4e+00 & 2.3e-06\\
      10 & 8.3e+04 & 5.2e+05 & 4.4e+02 & 3.2e+01 & 2.7e-06\\
      \hline
      12 & 5.2e+03 & 3.3e+04 & 2.9e+01 & 2.7e+00 & 1.9e-07\\
      12 & 2.1e+04 & 1.3e+05 & 1.2e+02 & 7.9e+00 & 3.1e-07\\
      12 & 8.3e+04 & 5.2e+05 & 4.7e+02 & 3.4e+01 & 3.5e-07\\
      \hline
    \end{tabular}
  \end{center}
  \caption{Numerical results of operator $S$ for the bean-shaped object.}
  \label{tbl:S2}
\end{table}


\begin{table}[ht!]
  \begin{center}
    \begin{tabular}{|c|cc|ccc|}
      \hline
      $m_c$ & $\omega$ & $n$ & $T_s$ & $T_a$ & $e$\\
      \hline
      6 & 5.3e+03 & 3.3e+04 & 2.2e+01 & 1.4e+00 & 3.3e-04\\
      6 & 2.1e+04 & 1.3e+05 & 8.6e+01 & 5.3e+00 & 5.1e-04\\
      6 & 8.5e+04 & 5.2e+05 & 3.5e+02 & 2.3e+01 & 4.7e-04\\
      \hline
      8 & 5.3e+03 & 3.3e+04 & 2.2e+01 & 1.9e+00 & 2.1e-05\\
      8 & 2.1e+04 & 1.3e+05 & 8.9e+01 & 5.8e+00 & 1.9e-05\\
      8 & 8.5e+04 & 5.2e+05 & 3.6e+02 & 2.4e+01 & 2.6e-05\\
      \hline
      10 & 5.3e+03 & 3.3e+04 & 2.2e+01 & 2.0e+00 & 6.8e-07\\
      10 & 2.1e+04 & 1.3e+05 & 9.1e+01 & 6.2e+00 & 4.8e-07\\
      10 & 8.5e+04 & 5.2e+05 & 3.7e+02 & 2.6e+01 & 8.3e-07\\
      \hline
      12 & 5.3e+03 & 3.3e+04 & 2.3e+01 & 2.2e+00 & 2.3e-08\\
      12 & 2.1e+04 & 1.3e+05 & 9.5e+01 & 6.7e+00 & 2.7e-08\\
      12 & 8.5e+04 & 5.2e+05 & 3.9e+02 & 2.9e+01 & 2.6e-08\\
      \hline
    \end{tabular}
  \end{center}
  \caption{Numerical results of operator $D$ for the ellipse.}
  \label{tbl:D1}
\end{table}

\begin{table}[ht!]
  \begin{center}
    \begin{tabular}{|c|cc|ccc|}
      \hline
      $m_c$ & $\omega$ & $n$ & $T_s$ & $T_a$ & $e$\\
      \hline
      6 & 5.2e+03 & 3.3e+04 & 2.8e+01 & 2.5e+00 & 1.2e-03\\
      6 & 2.1e+04 & 1.3e+05 & 1.1e+02 & 6.6e+00 & 2.0e-03\\
      6 & 8.3e+04 & 5.2e+05 & 4.5e+02 & 3.0e+01 & 1.0e-03\\
      \hline
      8 & 5.2e+03 & 3.3e+04 & 2.9e+01 & 2.3e+00 & 7.8e-05\\
      8 & 2.1e+04 & 1.3e+05 & 1.1e+02 & 7.0e+00 & 6.3e-05\\
      8 & 8.3e+04 & 5.2e+05 & 4.6e+02 & 3.0e+01 & 7.3e-05\\
      \hline
      10 & 5.2e+03 & 3.3e+04 & 3.0e+01 & 2.5e+00 & 1.6e-06\\
      10 & 2.1e+04 & 1.3e+05 & 1.2e+02 & 7.4e+00 & 4.1e-06\\
      10 & 8.3e+04 & 5.2e+05 & 4.7e+02 & 3.2e+01 & 4.0e-06\\
      \hline
      12 & 5.2e+03 & 3.3e+04 & 3.1e+01 & 2.7e+00 & 3.2e-07\\
      12 & 2.1e+04 & 1.3e+05 & 1.2e+02 & 8.0e+00 & 4.5e-07\\
      12 & 8.3e+04 & 5.2e+05 & 5.0e+02 & 3.4e+01 & 5.4e-07\\
      \hline
    \end{tabular}
  \end{center}
  \caption{Numerical results of operator $D$ for the bean-shaped object.}
  \label{tbl:D2}
\end{table}


\begin{table}[ht!]
  \begin{center}
    \begin{tabular}{|c|cc|ccc|}
      \hline
      $m_c$ & $\omega$ & $n$ & $T_s$ & $T_a$ & $e$\\
      \hline
      6 & 5.3e+03 & 3.3e+04 & 2.1e+01 & 2.3e+00 & 3.1e-04\\
      6 & 2.1e+04 & 1.3e+05 & 8.6e+01 & 5.4e+00 & 3.2e-04\\
      6 & 8.5e+04 & 5.2e+05 & 3.5e+02 & 2.3e+01 & 4.3e-04\\
      \hline
      8 & 5.3e+03 & 3.3e+04 & 2.1e+01 & 1.9e+00 & 1.5e-05\\
      8 & 2.1e+04 & 1.3e+05 & 8.8e+01 & 5.8e+00 & 1.7e-05\\
      8 & 8.5e+04 & 5.2e+05 & 3.5e+02 & 2.4e+01 & 1.0e-05\\
      \hline
      10 & 5.3e+03 & 3.3e+04 & 2.2e+01 & 2.0e+00 & 3.8e-07\\
      10 & 2.1e+04 & 1.3e+05 & 9.1e+01 & 6.2e+00 & 4.1e-07\\
      10 & 8.5e+04 & 5.2e+05 & 3.7e+02 & 2.6e+01 & 4.1e-07\\
      \hline
      12 & 5.3e+03 & 3.3e+04 & 2.3e+01 & 2.2e+00 & 1.6e-08\\
      12 & 2.1e+04 & 1.3e+05 & 9.6e+01 & 6.7e+00 & 3.3e-08\\
      12 & 8.5e+04 & 5.2e+05 & 3.9e+02 & 2.9e+01 & 1.8e-08\\
      \hline
    \end{tabular}
  \end{center}
  \caption{Numerical results of operator $D'$ for the ellipse.}
  \label{tbl:X1}
\end{table}

\begin{table}[ht!]
  \begin{center}
    \begin{tabular}{|c|cc|ccc|}
      \hline
      $m_c$ & $\omega$ & $n$ & $T_s$ & $T_a$ & $e$\\
      \hline
      6 & 5.2e+03 & 3.3e+04 & 2.8e+01 & 2.5e+00 & 1.6e-03\\
      6 & 2.1e+04 & 1.3e+05 & 1.1e+02 & 6.6e+00 & 1.5e-03\\
      6 & 8.3e+04 & 5.2e+05 & 4.5e+02 & 3.1e+01 & 7.6e-04\\
      \hline
      8 & 5.2e+03 & 3.3e+04 & 2.9e+01 & 2.3e+00 & 9.3e-05\\
      8 & 2.1e+04 & 1.3e+05 & 1.1e+02 & 7.0e+00 & 1.3e-04\\
      8 & 8.3e+04 & 5.2e+05 & 4.6e+02 & 3.0e+01 & 1.0e-04\\
      \hline
      10 & 5.2e+03 & 3.3e+04 & 3.0e+01 & 2.5e+00 & 6.3e-06\\
      10 & 2.1e+04 & 1.3e+05 & 1.2e+02 & 7.4e+00 & 5.0e-06\\
      10 & 8.3e+04 & 5.2e+05 & 4.7e+02 & 3.2e+01 & 3.6e-06\\
      \hline
      12 & 5.2e+03 & 3.3e+04 & 3.1e+01 & 2.7e+00 & 2.5e-07\\
      12 & 2.1e+04 & 1.3e+05 & 1.2e+02 & 7.9e+00 & 3.7e-07\\
      12 & 8.3e+04 & 5.2e+05 & 5.0e+02 & 3.5e+01 & 3.8e-07\\
      \hline
    \end{tabular}
  \end{center}
  \caption{Numerical results of operator $D'$ for the bean-shaped object.}
  \label{tbl:X2}
\end{table}


\begin{table}[ht!]
  \begin{center}
    \begin{tabular}{|c|cc|ccc|}
      \hline
      $m_c$ & $\omega$ & $n$ & $T_s$ & $T_a$ & $e$\\
      \hline
      6 & 5.3e+03 & 3.3e+04 & 2.4e+01 & 2.3e+00 & 6.9e-04\\
      6 & 2.1e+04 & 1.3e+05 & 1.0e+02 & 5.3e+00 & 4.2e-04\\
      6 & 8.5e+04 & 5.2e+05 & 4.0e+02 & 2.3e+01 & 3.1e-04\\
      \hline
      8 & 5.3e+03 & 3.3e+04 & 2.5e+01 & 1.9e+00 & 2.0e-05\\
      8 & 2.1e+04 & 1.3e+05 & 1.0e+02 & 5.8e+00 & 1.3e-05\\
      8 & 8.5e+04 & 5.2e+05 & 4.2e+02 & 2.4e+01 & 1.6e-05\\
      \hline
      10 & 5.3e+03 & 3.3e+04 & 2.7e+01 & 2.0e+00 & 4.5e-07\\
      10 & 2.1e+04 & 1.3e+05 & 1.1e+02 & 6.2e+00 & 1.3e-06\\
      10 & 8.5e+04 & 5.2e+05 & 4.4e+02 & 2.6e+01 & 7.4e-07\\
      \hline
      12 & 5.3e+03 & 3.3e+04 & 2.9e+01 & 2.2e+00 & 3.0e-08\\
      12 & 2.1e+04 & 1.3e+05 & 1.2e+02 & 6.7e+00 & 2.9e-08\\
      12 & 8.5e+04 & 5.2e+05 & 4.7e+02 & 2.8e+01 & 2.0e-08\\
      \hline
    \end{tabular}
  \end{center}
  \caption{Numerical results of operator $N$ for the ellipse.}
  \label{tbl:N1}
\end{table}

\begin{table}[ht!]
  \begin{center}
    \begin{tabular}{|c|cc|ccc|}
      \hline
      $m_c$ & $\omega$ & $n$ & $T_s$ & $T_a$ & $e$\\
      \hline
      6 & 5.2e+03 & 3.3e+04 & 3.2e+01 & 2.6e+00 & 1.2e-03\\
      6 & 2.1e+04 & 1.3e+05 & 1.3e+02 & 6.6e+00 & 1.1e-03\\
      6 & 8.3e+04 & 5.2e+05 & 5.1e+02 & 3.0e+01 & 1.4e-03\\
      \hline
      8 & 5.2e+03 & 3.3e+04 & 3.3e+01 & 2.4e+00 & 1.0e-04\\
      8 & 2.1e+04 & 1.3e+05 & 1.3e+02 & 7.0e+00 & 1.2e-04\\
      8 & 8.3e+04 & 5.2e+05 & 5.3e+02 & 3.0e+01 & 1.2e-04\\
      \hline
      10 & 5.2e+03 & 3.3e+04 & 3.5e+01 & 2.5e+00 & 4.2e-06\\
      10 & 2.1e+04 & 1.3e+05 & 1.4e+02 & 7.4e+00 & 4.6e-06\\
      10 & 8.3e+04 & 5.2e+05 & 5.6e+02 & 3.2e+01 & 3.4e-06\\
      \hline
      12 & 5.2e+03 & 3.3e+04 & 3.8e+01 & 2.8e+00 & 5.1e-07\\
      12 & 2.1e+04 & 1.3e+05 & 1.5e+02 & 8.0e+00 & 3.4e-07\\
      12 & 8.3e+04 & 5.2e+05 & 6.1e+02 & 3.4e+01 & 5.9e-07\\
      \hline
    \end{tabular}
  \end{center}
  \caption{Numerical results of operator $N$ for the bean-shaped object.}
  \label{tbl:N2}
\end{table}

The numerical results show that both the setup and the evaluation
algorithms scale linearly with respect to $\omega$ and $n$. The
log-square factor $\log^2\omega$ of the application algorithm is not
significant. The running time reported here are obtained with our
current Matlab implementation. We expect a careful C/C++ or Fortran
implementation to exhibit much lower absolute running times.  The
relative error is clearly controlled by the size of the Chebyshev grid
$m_c$. Even for relative small Chebyshev grids, the algorithms achieve
good accuracy. An attractive feature of the new algorithms is that
once the kernel function is provided the algorithms are fully
kernel-independent.

\section{Discussions}

This paper presented a new directional algorithm for rapid evaluation
of high frequency boundary integrals via oscillatory Chebyshev
interpolation and local FFTs. This algorithm is conceptually simple,
fast, and kernel-independent.

Among the previously presented algorithms, this algorithm shares with
the algorithm in \cite{engquist-2007-fdmaoc,engquist-2009-fdahfas} the
idea of using directional low-rank separated approximations. However,
in \cite{engquist-2007-fdmaoc,engquist-2009-fdahfas} the directional
upward and downward equivalent sources are computed in upward and
downward passes in a recursive fashion. In this algorithm, the
low-rank separated approximations are computed instead using local
FFTs combined with oscillatory Chebyshev interpolation. This makes
the current algorithm much simpler to implement.

Though the paper presented the algorithms in the case of smooth
boundaries, they work for piecewise smooth boundaries without much
modification. One assumption is that the discretization points are
equally spaced according to the arclength parameterization. This
assumption can be relaxed if the parameterization is a smooth function
of the arclength. For an irregular discretization, the non-uniform
FFTs \cite{beylkin-1995,dutt-1993,greengard-2004,potts-2001} can be
used instead. The main topic for future work is the extension of this
approach to the 3D boundary integral formulations of the scattering
problems.

\bibliography{ref.bib}

\end{document}